\documentclass[12pt]{article}
\newcommand{\Hb}{\hfil\break}
\newcommand{\Sumc}{\displaystyle{\sum_{(c)}}}
\newcommand{\Circ}{\bullet}
\newcommand{\Field}{{\bf k}}
\newcommand{\Jperp}{J^{\perp}}
\newcommand{\Hom}{{\rm Hom}}
\newcommand{\Real}{{\bf R}}
\newcommand{\End}{\nonumber\\}
\newcommand{\Lg}{{L}}
\newcommand{\Proof}{{\em Proof }}
\newcommand{\Ker}{{\rm ker}}
\newcommand{\Bbc}[1]{{\bf{#1}}}
\newcommand{\Comp}{\Bbc{C}}
\newcommand{\To}{\longrightarrow}
\title{Measuring comodules - their applications}
\author{Marjorie Batchelor\\
Department of Mathematics\\
King's College\\
Strand\\
London WC2R 2LS}
\date{June 1998\\}
\newtheorem{Definition}{Definition}[section]

\newtheorem{Lemma}[Definition]{Lemma}
\newtheorem{Proposition}[Definition]{Proposition}
\newtheorem{Examples}[Definition]{Examples}
\newtheorem{Remark}[Definition]{Remark}
\newtheorem{Remarks}[Definition]{Remarks}
\newtheorem{Corollary}[Definition]{Corollary}
\newtheorem{Construction}[Definition]{Construction}
\begin{document}
\maketitle
\begin{abstract}
Measuring comodules are defined and shown to provide a useful
generalization of the set of maps between modules with a broad
range of applications.  Three applications are described.
Connections on bundles are described in terms of measuring
comodules, enabling curvature to be defined under general algebraic
circumstances.  Loop algebras are realized via a short exact
sequence of measuring comodules, with the central extension given
by the curvature.  Finally dual comodules provide a method of
dualizing representations, which, when applied to representations
of loop algebras yield positive energy representations, and when
applied to representations of totally disconnected groups leads to
the smooth dual.
\end{abstract}
For some time measuring coalgebras have been employed as sets of
generalized maps between algebras [1].  The purpose of this paper
is to introduce measuring comodules which provide a set of
generalized module maps from  a module $M$ over an algebra $A$ to a
module $N$ over a different algebra $B$.

        The categorical implication of this construction are presented in [4].
This paper presents a more practical approach.  Not only is the construction
of categorical interest, but it has wide ranging potential applications.  Three
are described here.

        The first describes connections on bundles.  A connection is that
construction which is required to describe covariant  differentiation of a
section of a vector bundle.  As such it is amenable to algebraic description,
and indeed gives an example of a measuring comodule.  The curvature of a
connection can likewise be described as an element of a measuring comodule.

        The second application generalizes an alternative construction of the
universal enveloping algebra of a Lie algebra using measuring
coalgebras. That construction proceeds as follows.  Given a
representation of a Lie algebra $L$ as derivations of an algebra
$A$, the universal enveloping algebra $UL$ arises as the
subcoalgebra-subalgebra of the universal measuring coalgebra
$P(A,A$) generated by $L$.  In an earlier paper quantum group -
like objects were shown to arise by considering
subcalgebra-subalgebras of $P(A,A)$ generated by certain sets of
difference operators [2]. Here similar sets of difference operators
are used to generate subcomodule-subalgebras of the universal
measuring comodule $Q(M,M)$ for a suitable $A$ module $M$.  The
resulting algebras are closely related to loop algebras and their
central extensions. The cocycle defining the central extension
arises as the trace of a curvature.

        The last application concerns the dual comodule of a $A$ module $M$,
$Q(M,\Comp)$, where $A$ is an algebra over $\Comp$.  This comodule
itself becomes an $A$ module with a strong finiteness property:
every element is contained in a finite dimensional $A$ submodule.
Two examples are considered.    In the first case $M$ is a
representation of a totally disconnected group $G$, for example, an
algebraic group over the $p$-adic numbers.  The construction of
interest concerns the dual comodule of $M$ considered as a module
for the group algebra $A = \Comp K$ where $K$ is a compact open
subgroup of $G$.  The resulting dual comodule turns out not only to
be a representation of $K$, but of the whole of $G$.  It is closely
related to the smooth dual.

        In the second example $M$ is a level $k$ representation of a loop algebra
$\Lg[x,x^{-1}]$ (where $L$ is a finite dimensional semi-simple Lie
algebra). The dual comodule for $M$ considered as an $A = U\Lg[x]$
module is not only an $\Lg[x]$ module but a level k representation
of $\Lg[x,x^{-1}]$. Moreover it contains, as a functorially
identifiable submodule a positive energy piece.  Where $M$ is
positive energy, this piece is the dual positive energy level $k$
representation of $\Lg[x,x^{-1}]$.

        The paper is organised as follows.  Section 1 describes measuring
coalgebras and comodules.  While essentially the constructions are the same
as those described in [4], here they are presented in a simplified algebraic
context rather than the more general categorical setting.  Connections are
described in section 2,  the association with loop algebras in section 3, and
finally the two applications of dual comodules in section 4.

        I would like particularly to thank Martin Hyland.  Much of the work
described here grew out of regular conversations held during my period of
retirement.  I would also like to thank Shaun Stevens for lessons in
disconnected groups.  Finally, I owe a great debt of gratitude to Alice Rogers
and King's College, London who have funded my return to active
mathematics, and caused this work to reach printed form.

\section{Measuring coalgebras and measuring comodules}

        Although measuring coalgebras (and dual coalgebras in particular)
have been around for a long time, I will develop the theory of
measuring coalgebras and measuring comodules in parallel, as the
first serves as an accessible model for the second.

\begin{Definition}{\rm Measuring coalgebras.  If $A$ and $B$ are algebras over a
field $\Field$ a {\em measuring coalgebra}  is a coalgebra $C$ over
$\Field$ with comultiplication
\begin{equation}
  \Delta : C \to C \otimes C, \quad \Delta c = \Sumc c_{(2)} \otimes c_{(1)}
\end{equation}

and counit $\epsilon: C \to \Field$ together with a linear map,
called a {\em measuring map}
  \begin{equation}
    f: C \to \Hom_{\Field}(A,B)
  \end{equation}

such that
 \begin{enumerate}
 \renewcommand{\labelenumi}{(\roman{enumi})}
  \item $\phi c(aa')  =   \Sumc \phi c_{(2)}(a)\phi c_{(1)}(a')$
  \item $\phi c(1_A) =  \epsilon(c)1_B$
 \end{enumerate}
for $a, a'$ in $A$, $1_A, 1_B$, the appropriate identity elements.
The map $\phi$ is said to measure.

Statements i and ii are equivalent to the statement that the
transpose map
 \begin{equation}
                \phi: A \to \Hom_{\Field}(C, B)
  \end{equation}
is an algebra homomorphism where the multiplication in
$\Hom_{\Field}(C,B)$ is given by
 \begin{equation}
                \mu\Circ \nu (c)  =   \Sumc\mu(c_{(2)})\nu(c_{(1)})
 \end{equation}

with identity
\begin{equation}
  1(c)  =  \epsilon(c)1_B.
\end{equation}

        The following proposition summarizes results about measuring
coalgebras described in [4].
 }\end{Definition}

  \begin{Proposition} \Hb {\rm
 \begin{enumerate}
 \renewcommand{\labelenumi}{(\roman{enumi})}
\item Given algebras $A$, $B$, there is a category of measuring
coalgebras $C(A,B)$ whose objects are measuring coalgebras
$(C,\phi)$ and whose maps \hfil $r: (C,\phi) \to (C',\phi')$ are
 coalgebra maps $r:C \to C'$ such that the diagram
\setlength{\unitlength}{0.7cm}
\begin{picture}(12,5)
 \put(3,3.5){$C$}
 \put(4,3.5){\vector(1,0){5}}
 \put(10,3.5){$\Hom_{\Field}(A,B)$}
 \put(4,3){\vector(1,-1){1.8}}
 \put(8.0,1.0){\vector(1,1){1.8}}
 \put(6.5,0.5){$C'$}
 \end{picture}

commutes.

        \item  The subcategory of finite dimensional measuring coalgebras is
dense in $C(A,B)$.  Essentially, every measuring coalgebra is a
limit of finite dimensional subcoalgebras.  (For a discussion of
density see [6], chapter 5.)

        \item  The category $C(A,B)$ has a final object, $(P(A,B), \pi)$ called the
universal measuring coalgebra.

Thus there is a correspondences of sets
 \begin{equation}
 \mbox{Coalgebra maps} (C, P(A,B))  \longleftrightarrow \mbox{ Algebra maps} (A,\Hom_{\Field}(C,B)).
 \end{equation}

        \item  If $A_i, i, = 1, 2, 3$ are algebras there is a map

               $ m:C(A_2,A_3)\times C(A_1,A_2) \to C(A_1,A_3)$.

In particular, $P(A,A)$ is a bialgebra.

        \item The universal measuring coalgebra $P(A,A)$ is a bialgebra.
        \end{enumerate}

Proofs.  Full proofs can be found in [4]  However, as the presentation in [4] is
highly categorical and more general than is necessary here, direct proofs of ii
and iii are indicated here.

        ii.To establish density  it is sufficient to show that every element $c$ of a
coalgebra $C$ is contained in a finite dimensional subcoalgebra
$C_1$. As this result is the essential property of coalgebras, the
proof, from [9] p46  is repeated here.
        Let $C' = \Hom_{\Field}(C,k)$ be the dual algebra and consider the action of $C'$
on $C$ given by
\begin{equation}
   c'\cdot c  =  \Sumc c'(c_{(2)})c_{(1)}
\end{equation}

Evidently the $C'$ module $V$ generated by $c$ is finite
dimensional, and
\begin{equation}
   C' \to \mbox{End}(V)
\end{equation}
is an algebra homomorphism of cofinite dimensional kernel $J$.

        Let $\Jperp$ be the subspace of $C$ on which $J$ is identically zero.  Finally
notice that $\Jperp \to \Hom(C'/J, \Field)$ is finite dimensional
and $c$ is in $\Jperp$.  A subcoalgebra of a measuring coalgebra is
itself a measuring coalgebra(with the restriction of the measuring
map), hence the result.

        iii.  This depends on two categorical properties of coalgebras.

        a)  Arbitrary coproducts exist in the category of coalgebras and
coalgebra maps.

        b)  Coequalizers also exist in this category.

The construction of $P(A,B)$ proceeds as follows.  Consider the
collection $\{(C_{\lambda},\phi_{\lambda})\}$ of finite dimensional
measuring coalgebras. Form the coproduct $\sqcup_{\lambda}
C\lambda$.  This is a measuring coalgebra. Now consider the set
$\{\rho({\lambda},\mu)\}$ of maps $\rho({\lambda},\mu):C_{\lambda}
\to C_{\mu}$ of finite dimensional measuring coalgebras. Form
$\sqcup_{\rho(\lambda,\mu)}C_{\lambda}$.  This is also a measuring
coalgebra.  There are two maps

\begin{equation}
  \alpha,\beta:\sqcup_{\rho(\lambda,\mu)}   C_{\lambda}    \to   \sqcup_{\lambda}C_{\lambda}.
\end{equation}

 On $C_{\lambda}$, $\alpha$ is just the inclusion
 $C_{\lambda} \to \sqcup{}_{\lambda}C_{\lambda}$
 while $\beta$ is the composition of
 $\rho(\lambda,\mu)$ with the inclusion $C_{\mu} \to \sqcup{}_{\lambda}C_{\lambda}$.

        The claim is that the coequalizer $P(A,B)$ has the desired universal
property.  If $(D,\psi)$ is a measuring coalgebra then $D$ is the
union of finite dimensional subcoalgebras $D_{\nu}$.  Evidently
there is a map $r_{\nu}:D_{\nu} \to P(A,B)$ and this map is unique.
The uniqueness of $r_{\nu}$ guarantees that the map $\rho:D
\to P(A,B)$ given by
 $\rho(d) = r_{\nu}(d)$ if $d$ is in $D_{\nu}$ is well defined.
 }\end{Proposition}

 \begin{Examples}\label{EXAMP}\Hb{\rm
 \begin{enumerate}
 \renewcommand{\labelenumi}{(\roman{enumi})}
 \item $P(A,B)$ is intended to generalize the set of
algebra homomorphisms from $A$ to $B$, and so it does.  Let
 $C_0  = \Field g$ be the one dimensional coalgebra with
 $\Delta g = g\otimes  g$, $\epsilon(g) = 1$.
Then a map $\phi: C_0 \to \Hom_{\Field}(A,B)$ measures if and only
if $\phi(g)$ is an algebra homomorphism.  Thus $P(A,B)$ contains
all algebra homomorphisms.

 \item        Let $C_1 = \Field g \oplus \Field \gamma$, $g$ as above, and
 let $\Delta \gamma = g \otimes \gamma + \gamma \otimes g$,  $\epsilon(\gamma) = 0$.
Then $\phi: C_1 \to \Hom_{\Field}(A,B)$ measures if and only if
$\phi(g)$ is an algebra homomorhpism and $\phi(\gamma)$ is a
derivation with respect to $\phi(g)$. That is,

\begin{equation}
   \phi(\gamma)(aa')  =  \phi(\gamma)(a)\phi(g)(a') + \phi(g)(a)\phi(\gamma)(a').
\end{equation}

 \item          More generally if $L$ is a Lie algebra over $\Field$, then $L\oplus C_0$
  can be
  given the structure of coalgebra with comultiplication
  $\Delta \gamma = \gamma \otimes g + g \otimes \gamma$ and
$\epsilon(\gamma) = 0$ for $\gamma$ in $L$.  Suppose
\begin{equation}
 \phi: L\oplus C_0 \to \Hom_{\Field}(A,A)
\end{equation}
is a measuring map such that
\begin{equation}
  \phi[\nu,\gamma]  =  [\phi\nu, \phi\gamma],  \qquad   \phi(g)  =  Id
\end{equation}
for $\nu, \gamma$ in $L$.  By the universal property of $P(A,A)$
there is a map of measuring coalgebras
\begin{equation}
     \rho: L \oplus C_0 \to P(A,A).
\end{equation}
However $P(A,A)$ is an algebra, and in fact the following is true.
\end{enumerate}
}
\end{Examples}
\begin{Proposition}\label{ULprop}  {\rm If the map $\phi$ is injective on $L$ the subalgebra
of $P(A,A)$ generated by the image of $\rho$ is isomorphic to the
universal enveloping algebra $UL$.

        The proof follows from the universal property of $P(A,A)$ and facts
about bialgebras [2].  In that paper I considered subalgebras of
$P(A,A)$ generated by measuring coalgebras $L \oplus \Comp K$,
where $K$ is a group, $\Comp K$ has the usual comultiplication
 $\Delta  k = k \otimes k$ for $k$ in $K$, and elements of $L$ have a slightly skew
version of the usual comultiplication for derivations,
\begin{equation}
   \Delta E = E \otimes k + k^{-1} \otimes E
\end{equation}
which is characteristic of difference operators.  These objects
resemble quantum groups.  The construction in section \ref{GENsec}
of this paper uses the same procedure to construct subalgebras of
the universal measuring comodule (defined below) which are related
to central extensions of loop algebras.
 }\end{Proposition}
 \begin{Definition}{\rm {\em Measuring comodules.}
 Let $M$ be an $A$ module and let $N$ be a $B$ module (all
modules and algebras are vector spaces over $\Field$).  When it is
necessary to emphasize the algebra over which $M$ and $N$ are
modules write ${}^AM$, ${}^BN$.  Let $(C,\phi)$ be a measuring
coalgebra in $C(A,B)$. Recall that a {\em comodule} over $C$ is a
vector space with a comultiplication
\begin{equation}
  \Delta: D \to C \otimes D,       \quad   \Delta(d) = \sum_{(d)} d_{(1)} \otimes  d_{(0)}
\end{equation}
In addition I will assume that $(\epsilon \otimes1)\Delta = 1$.
When it is necessary to keep track of the coalgebra over which  $D$
is a comodule, write ${}_CD$.  A $\Field$-linear map
 \begin{equation}
     \psi:D \to \Hom_{\Field}(M,N)
 \end{equation}
{\em measures} if
\begin{equation}
     \psi(am)  =  \sum_{(d)}\phi d_{(1)}(a) \psi d_{(0)}(m).
\end{equation}
The pair $(D,\psi)$ is called a measuring comodule, and $\psi$ is
called a measuring map.

        Equivalently $\psi$ measures if and only if the corresponding  transpose
map
\begin{equation}
                        \psi: M \to \Hom_{\Field}(D,N)
\end{equation}
is a map of $A$ modules, where the $A$ module structure on
$\Hom_{\Field}(D,N)$ is given by
\begin{equation}
                a\Circ \beta(d)  =  \sum_{(d)} \phi d_{(1)}(a) \beta d_{(0)}.
\end{equation}

        Again, results from [4] are summarized in the following proposition.
 }\end{Definition}
 \begin{Proposition} \Hb{\rm
 \begin{enumerate}
   \renewcommand{\labelenumi}{(\roman{enumi})}
  \item Given a measuring coalgebra  $C$ in $C(A,B)$
there is a category ${}_CD(M,N)$ whose objects are measuring
comodules $(D,\psi)$ and whose maps $\sigma:(D,\psi) \to
(D',\psi')$ are comodule maps
 $\sigma: D \to D'$ such that  the diagram
\par
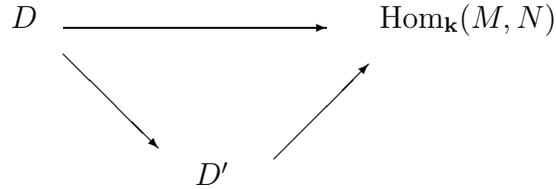
\setlength{\unitlength}{0.7cm}
 \begin{picture}(12,5)
 \put(3,3.5){$D$}
 \put(4,3.5){\vector(1,0){5}}
 \put(10,3.5){$\Hom_{\Field}(M,N)$}
 \put(4,3){\vector(1,-1){1.8}}
 \put(8.0,1.0){\vector(1,1){1.8}}
 \put(6.5,0.5){$D'$}
 \end{picture}
\par
commutes.

 \item  The subcategory of ${}_CD(M,N)$ whose objects are the finite
dimensional measuring comodules is a dense subcateory of
${}_CD(M,N)$.

 \item The category ${}_CD(M,N)$ has a final object, ${}_CQ(M,N)$.  This has the
property that there is a correspondence
\begin{equation}
        C-\mbox{comodule maps}(D, {}_CQ(M,N)) \longleftrightarrow A-\mbox{module maps}(M, \Hom(D,N)).
\end{equation}
 \item If $M_i$  are modules over algebras $A_i, i = 1,2,3$, and if $C$, $C'$ are in
$C(A_1,A_2), C(A_2,A_3)$ respectively, then there is a map
\begin{equation}
                 {}_CD(M_2,M_3)\times {}_{C'}D(M_1,M_2) \stackrel{m}{\to} _{m(C \times C')}D(M_1,M_3).
\end{equation}
In particular, ${}_CQ(M,M)$ is a comodule algebra, for $C \to
\Hom(A,A)$ a measuring coalgebra, $M$ an $A$ module.
\item  If $A = B$ and $M = N$, then ${}_CQ(M,N)$ is a comodule algebra.
 \end{enumerate}

Proofs. ii.   Again the important step is to show that if $D$ is a
$C$ comodule then each $d$ in $D$ is contained in a finite
dimensional subcomodule.

        Define an action $\Circ$ of $C'$, the linear dual of $C$ on $D$ via
\begin{equation}
                        a\Circ d   =  \sum_{(d)} a(d_{(1)})d_{(0)}.
\end{equation}
Choose an element $d$ of $d_0$ and let $D_0 = C'\Circ d_0$.
Evidently $D_0$ is a finite dimensional $C'$ module and $d_0 =
1\Circ d_0$ is in $D_0$.

        The full linear dual $D'$ (of $D$) is also a $C'$ module, with the action given
explicitly by
\begin{equation}
                        a*d(d)  =  \sum_{(d)} a(d_{(1)})d(d_{(0)}).
\end{equation}
The subset $D_0^{\perp}  =  \{d \epsilon D' : d(D_0) = 0\}$ is a
submodule, and $D_0' = D'/D_0{}^{\perp}$.

        But $(D_0^{\perp} )^{\perp}$   is then a subcomodule of $D$, and $D$ is contained in
        $(D_0{}^{\perp})^{\perp}$.  Since  $(D_0{}^{\perp})^{\perp}$  includes in $(D'/D_0^{\perp})'$,
$(D_0^{\perp})^{\perp}$ must be a finite dimensional comodule as
required.
        The proofs of i and iii are identical in format to the corresponding
statements for measuring coalgebras.
 }\end{Proposition}
\begin{Examples}\Hb{\rm
\begin{enumerate}
  \renewcommand{\labelenumi}{(\roman{enumi})}
 \item  Let $C_0 = \Field g$ as in example \ref{EXAMP}.i and suppose that
$\phi:C \to \Hom_{\Field}(A,B)$ measures, so that $\phi(g)$ is an
algebra homomorphism. Let $D$ be the comodule with
 $D  = \Field d$ and comultiplication $\Delta d = g \otimes d$. Let
  $\psi:D \to \Hom_{\Field}(M,N)$
be a linear map.
        Recall that the pullback of $N$, $\phi(g)^*N$ is an $A$ module. Then $\psi$ measures
if and only if
\begin{equation}
                        \psi(d): M \to \phi(g)^{*}N
\end{equation}
is a map of $A$ modules.

 \item  If $A = B$ and if $C$ contains the measuring comodule $C_0$ with
  $\phi(g) = 1$, the measuring comodule ${}_CQ(M,N)$ contains the vector space $H$ of
all genuine $A$ module maps from $M$ to $N$ as follows.

        Any vector space, for example $H$, is trivially a $C$ comodule with
comultiplication $\Delta h = g \otimes h$.  The inclusion
 $\psi: H \to \Hom_{\Field}(M,N)$ is then a measuring map.  By the universal
property there is a unique map of measuring comodules
 $\rho: H \to {}_CQ(M,N)$.  Since $\psi$ is an inclusion, so must $\rho$ be.

  \item  Any algebra  can be considered as a module over itself acting by
left multiplication.  If $C \to \Hom(A,B)$ is a measuring
coalgebra, by considering $C$ as a comodule over itself,
 $C \to \Hom(A,B)$ is also a measuring comodule.

 \item  For an element $a$ of an algebra $A$ let $\iota_a$ denote the  inner derivation
\begin{equation}
                \iota_a(b)  = [a,b]  = ab - ba
\end{equation}
Let $I_A$ denote the Lie algebra of inner derivations of $A$.  As
in example 1.3.iii, let $C$ be the measuring coalgebra
  $C = I_A \oplus C_0, C  \to \Hom_{\Field}(A,A)$. Now put a $C$ comodule structure on $A$,
\begin{equation}
                \Delta(a) = g \otimes a  +  \iota_a \otimes 1.
\end{equation}
Now let $M$ be an $A$ module.  With the comodule structure above
the inclusion $A \to \Hom_{\Field}(M,M)$ sending $a$ to left
multiplication by $a$ gives $A$ the structure of a measuring
comodule. This construction generalizes the observation that for
modules over commutative rings, left multiplication is a module
map.
 \end{enumerate}
    }
\end{Examples}
\begin{Remarks}\Hb{\rm
 \begin{enumerate}
  \renewcommand{\labelenumi}{(\roman{enumi})}
 \item If $\tau:(C,f) \to (C',f')$ is a map of measuring coalgebras,
in particular $\tau$ is a comodule map, so that ${}_CQ(M,N)$ can be
considered as a $C'$ comodule.  Since $\tau$ is a map of measuring
coalgebras ${}_CQ(M,N)$ is in fact in ${}_{C'}D (M,N)$, and hence
by the universal property there is a unique map
 ${}_CQ(M,N) \to {}_{C'}Q(M,N)$. All universal measuring comodules ${}_CQ(M,N)$ thus map to
$_{P(A,B)}Q(M,N)$, which will often be denoted $Q(M,N)$.

\item  The construction $Q(M,N)$ serves as the set of ``module maps from
an $A$ module $M$ to a $B$ module $N$'' even when $A$ is not the
same as $B$. The paper [4] arose from the desire to put this
curiosity into a sound categorical context.
 \end{enumerate}
 }\end{Remarks}
\section{Connections}

        Given an algebra $A$ of functions, and a set $V$ of derivations of $A$(vector
fields), a connection is that which is needed to define covariant
differentiation by elements of $V$ on a module $M$ (for example,
sections of a bundle) over $A$.
        This is a completely algebraic statement and as such lends itself to
restatement in terms of measuring comodules.

\begin{Definition}{\rm
 {\em Loose connections} Let $A$ be an algebra and let $M$ be a module over $A$. Let $C$
be a measuring coalgebra, and let D be a comodule over $C$ which is
also an $A$ module.  A loose connection is a measuring map
\begin{equation}
                        \nabla: D \to \Hom_{\Field}(M,M)
\end{equation}
which additionally satisfies the requirement that $\nabla$ be a map
of $A$ modules in the sense that
\begin{equation}
                        \nabla(a\xi)(m) = a\nabla\xi(m).
\end{equation}
  }\end{Definition}
 \begin{Examples}\Hb{\rm
  \begin{enumerate}
    \renewcommand{\labelenumi}{(\roman{enumi})}
  \item Connections on a vector bundle.  Let
 $A = C^{\infty}(Y)$ where $Y$ is a smooth manifold and Let $V$ be the Lie algebra
of vector fields on $Y$ and let $C = V \oplus \Comp 1$.  Let
 $D = V \oplus A$ with the comultiplication
\begin{eqnarray}
        \Delta&:& D \to C \otimes D, \End
        \Delta(\psi) &=& 1\otimes \psi + \psi\otimes 1, \psi \in V\End
         \Delta(a)& =& 1\otimes a  + \iota_a \otimes 1 a \in A
\end{eqnarray}
Notice that $D$ is an $A$ module.  Let $E$ be a vector bundle over
$Y$ and let $\Gamma(Y,E)$ denote the smooth sections of $E$ over
 $Y$. Thus $M = \Gamma(Y,E)$ is a module for $A$.  In this setting loose connections
are precisely Koszul connections (see [8]).

 \item  Connections on a principle bundle. (See [7].) Let $Y$ be a manifold
and let $P$ be a principle $G$ bundle over $Y$.  Let
 $M = C^{\infty}(P)$. Observe that $C^{\infty}(Y)$ includes in $M$ as those
functions which are constant on the fibres of $P$, hence $M$ is a
 $C^{\infty}(Y)$ module. In addition the group algebra $\Real G$ acts on $M$ via
right translation. The action of $\Real G$ commutes with the action
of $C^{\infty}(Y)$.

        Let $A  =  C^{\infty}(Y) \otimes\Real  G$.  Let $V$ be the Lie algebra of vector fields on $Y$.
Observe that the coalgebra $C$ above becomes a measuring coalgebra
with measuring map
\begin{eqnarray}
        \phi: C &\To &\Hom \left(C^{\infty}(Y) \otimes\Real  G,C^{\infty}(Y) \otimes\Real  G \right), \End
        \phi(\psi)f \otimes  g &=& \psi f \otimes g.
        \end{eqnarray}
Let $D$ be the comodule of example i above.  This is a
$C^{\infty}(Y)\otimes \Real  G$ module with the trivial action of
$G$ on $C^{\infty}(Y)$ and $V$. Also notice that $D$ contains $C$
as a subcomodule (with its usual coproduct).  A loose connection
\begin{equation}
                \nabla: D \to \Hom(C^{\infty}(P),C^{\infty}(P))
\end{equation}
in this setting corresponds to a connection on the principle bundle
if and only if additionally $\nabla$ restricted to the subspace $C$
defines a measuring coalgebra
\begin{equation}
                \nabla: C \to \Hom(C^{\infty} (P),C^{\infty}(P)).
\end{equation}
 \end{enumerate}
 }
 \end{Examples}
\subsection{Curvature} {\rm
Curvature can be defined for any measuring comodule equipped with
Lie bracket, in particular for loose connections. Recall that any
measuring comodule, $D$ in particular, comes with a map of
measuring comodules
\begin{equation}
                        \rho:D \to Q(M,M).
\end{equation}
Recall that Q(M,M) is a comodule-algebra.  If $D$ contains a
subspace $V$ on which a Lie bracket is given, for $\xi,\psi$ in
 $V$, write
 \begin{equation}
               \Omega(\xi,\psi) = \rho(\xi)\rho(\psi) - \rho(\psi)\rho(\xi) - \rho([\xi,\psi]).
 \end{equation}
This map
 \begin{equation}
                \Omega: V\otimes V \to Q(M,M)
 \end{equation}
is the {\em curvature} of the loose connection $\nabla$ on V.

 \begin{Remark} {\rm   In all classical cases, the coalgebra $C$ is always
  $V \oplus \Comp 1$, the comodule $D$ is always $V \oplus A$, where $V$ is the Lie algebra of
derivations of $A$, and one is only interested in the restriction
of $\nabla$ to $V$.  There is no harm, however, in allowing this
greater generality.  In the next section a very different example
demonstrates the advantages of being broad minded.

        This section concludes with a result which is well known for
conventional connections.
 }\end{Remark}
  \begin{Proposition} {\rm   If V is a set of primitive elements, ie, with comultiplication
  \begin{equation}
                \Delta x = \xi \otimes 1 + 1 \otimes \xi
  \end{equation}
then $\Omega(\xi,\psi)$ determines a module map
  \begin{equation}
                \Omega(\xi,\psi): M \to M.
  \end{equation}
  }\end{Proposition}
\Proof.  Direct calculation (observing that Q(M,M) is a comodule algebra, ie,
that multiplication preserves the comodule structure) shows that
  \begin{equation}
                \Delta \Omega(\xi,\psi) = 1\otimes \Omega(\xi,\psi).
  \end{equation}
But this is exactly the statement that $\Omega(\xi,\psi)$ is a
module map.

\section{Generalizations of universal enveloping algebras using measuring comodules}
\label{GENsec}

        In proposition \ref{ULprop} the universal enveloping algebra is constructed as a
subalgebra of the universal measuring coalgebra generated by a Lie algebra of
derivations.  The original Lie algebra can be identified as the subspace of
primitive elements.

        This construction can be generalized, replacing primitive elements
(derivations) with elements $E$ of a measuring coalgebra with the
assymetric comultiplication $\Delta E = E \otimes K + K^{-1}
\otimes E$, where $K$ is an (invertible) group like element.  The
algebras generated by such $E$ and $K$ resemble quantum groups, and
were the subject of [2].

        This construction is now generalized again, replacing the universal
measuring coalgebra with the universal measuring comodule.  The
resulting algebra has, as the analogue of its Lie algebra of
primitive elements, a Lie algebra related to the central extensions
of loop algebras.  The central term arises as the trace of the
curvature.

\begin{Construction}{\rm
Let $A$ be an algebra and let $M$ be an $A$ module.  Let $P_0$ be a
subcoalgebra subalgebra of $P(A,A)$, and define
\begin{eqnarray}
                V_0 &=& \{v \epsilon Q(M,M): \Delta v \epsilon P_0 \otimes Q(M,M)\}. \End
        V &=& \{v \epsilon Q(M,M), \Delta v \epsilon P(A,A)\otimes 1 + 1
        \otimes Q(M,M) + P_0 \otimes Q(M,M)\}.\End
\end{eqnarray}
Evidently $A$ is contained in $V$.   The subcomodule $V$ is the
generalization of primitive elements referred to above.

 Suppose a ``trace''
 \begin{equation}
                \tau: V_0 \to \Comp
\end{equation}
is given. Let $K$ be the kernel of $\tau$.  Define
\begin{eqnarray}
                V_{0\tau} &=& \{v \epsilon V_0: [K,v] \leq K\} \End
                V_{\tau} &=& \{v \epsilon V: [V_{0\tau},v] \leq K\}
\end{eqnarray}
Observe that $V_{\tau}$ is not an algebra, but the Jacobi identity
guarantees that it is closed under Lie bracket.  $V_0$ is a
subalgebra of $Q(M,M)$ and hence $V_{0\tau}$ is a Lie subalgebra.
It is not hard to check that there is a short exact sequence of Lie
algebras
\begin{equation}
                0 \to V_{0\tau}/K \to V_{\tau}/K \to V_{\tau}/V_{0\tau} \to 0.
\end{equation}
Moreover $V_{\tau}/K$ is a central extension of
$V_{\tau}/V_{0\tau}$.

        Suppose now $\mu:V_{\tau}/V_{0\tau} \to V_{\tau}$ is any linear section of the
projection $V_{\tau} \to V_{\tau}/V_{0\tau}$. The image
$\mu(V_{\tau}/V_{0\tau})$ inherits a Lie bracket from
$V_{\tau}/V_{0\tau}$: hence the associated curvature $\Omega_{\mu}$
takes values in $V_{0\tau}$. While the $\Omega_{\mu}$ may depend on
the section $\mu$, the trace of the curvature does not. In fact
$V_{\tau}/K$ is the central extension of $V_{\tau}/V_{0\tau}$ with
cocycle $c$ defined by
\begin{equation}
                 c(v,w) = \tau(\Omega_{\mu}(\mu v,\mu w))
\end{equation}
for $v, w$ in $V_{\tau}/V_{0}$. Familiar  examples arise from
looking at particular subspaces of $V$.
 }
 \end{Construction}
\begin{Examples}\Hb{\rm
 \begin{enumerate}
 \renewcommand{\labelenumi}{(\roman{enumi})}
  \item   Let $A = M = \Comp[x]$.    Let
\begin{equation}
                P_0  = \{p \epsilon P(\Comp[x],C): p(x^n) = 0 \mbox{ for almost all } n\}.
\end{equation}
Explicitly $P_0$ has basis $\{\beta_{j}\}$ with comultiplication
and measuring map given by
\begin{equation}
                \Delta \beta_{j}
                = \sum_{k=0} \beta_k \otimes\beta_{j-k},
                \qquad \phi(\beta_{j}) = j! \frac{d}{dx}|_0  .
\end{equation}
Define a trace $\tau$ on $V_0$
\begin{equation}
                \tau(v) = \sum_{j=0}^{\infty} \beta_{j}(v(x^j)).
\end{equation}
To see this is well defined, observe that for $v$ in $V_0$
\begin{equation}
                v(x^n) = \sum_{(v)}v_{(1)}(x^n)v_{(0)}(1).
\end{equation}
Thus $\beta_{j}(v(x^n)) = 0 $ for greater than the greatest degree
of the $v_{(0)}(1)$, and the sum defining $\tau$ is always a finite
sum.

        This example contains two very well known examples as Lie
subalgebras.  First consider $\Comp[\alpha^{-1}]$.  This is can be
given the structure of a coalgebra with
\begin{eqnarray}
                \Delta(\alpha^{-i}) &=& \alpha^{-i}\otimes 1
                + \sum_{k=0}^{i-1} \beta_k \otimes \alpha^{-i+k},  i>0. \End
                \epsilon(\alpha^{i}) &=& \delta_{i,0}.
\end{eqnarray}
Define a map $\phi:\Comp[\alpha^{-1}] \to \Hom(\Comp[x],\Comp[x])$
via
 \begin{equation}
                \phi(\alpha^{i})x^n =  \left\{ \begin{array}
                  {r@{\quad}l}
                  x^{i+n}    & i + n \geq 0 \\
                   0 & \mbox{otherwise}.
                                                 \end{array} \right.
 \end{equation}
It is routine, if surprising, to verify that this map measures. Now
$\Comp[\alpha,\alpha^{-1}]$ can be considered a comodule over
$\Comp[\alpha^{-1}] \oplus P_0$ with comultiplication given by
\begin{equation}
                \Delta \alpha^{i} = 1 \otimes \alpha^{i}
\end{equation}
for $i>0$ otherwise as above.  Clearly the measuring map $\phi$
above extends to all of $\Comp[\alpha,\alpha^{-1}]$.  It is not
hard to check that the image of $\Comp[\alpha,\alpha^{-1}]$ lies in
$V_{\tau}$. The image of $\Comp[\alpha,\alpha^{-1}]$ in
$V_{\tau}/K$ is the familiar central extension of the abelian Lie
algebra $\Comp[\alpha,\alpha^{-1}]$ with cocycle
\begin{equation}
                c(\alpha^k,\alpha^j) = k\delta_{k,-j}.
\end{equation}
 \item  With $P_0$ and $\tau$ as before let $T$ be the vector
space with basis $\{T_i, i \in \Bbc{Z}\}$, and put a comodule
structure on $T$ via
\begin{equation}
                \Delta(T_i) = T_i \otimes 1 + 1 \otimes T_i +
                \sum_{i+k<0} k \beta_k \otimes \alpha^{k+i} + \beta^k \otimes T_{k+i} .
\end{equation}
Observe that $T \oplus P_0\oplus \Comp[\alpha^{-1} ]$ is in fact a
coalgebra if the counit on $T$ is defined to be identically $0$.
Extending $\phi$ of the previous example via
\begin{equation}
               \phi(T_i)(x^n)  =  \left\{ \begin{array}
                  {r@{\quad}l}
                  x^{i+n}    & i + n \geq 0 \\
                   0 & \mbox{otherwise}.
                                                 \end{array} \right.
 \end{equation}

gives $T \oplus P_0\oplus \Comp[\alpha^{-1} ]$  the structure of a
measuring coalgebra, and hence a measuring comodule.  Again the
image lies in $V_{\tau}$.   The image of $T$ in
$V_{\tau}/V_{0\tau}$ is isomorphic to the Lie algebra of
derivations of $\Comp[x,x^{-1}]$, and its image in $V_{\tau}/K$ is
the variant of the Virasoro algebra with cocycle
\begin{equation}
                c(T_m,T_n) =  \frac16 (m^3-m) \delta_{m,-n}.
\end{equation}
 \item  Return now to a general algebra $A$, and suppose
that $M = A$, and suppose also that $\tau$ is given so that $K$,
$V_{\tau}$ and $V_{0\tau}$ are as described in 3.1. Let $\Lg$ be a
finite dimensional (semi-simple) Lie algebra, which is faithfully
represented by $\rho: \Lg \to \mbox{End}(W)$.  Let $M(W) = A
\otimes W$. Then the identification of
 $\Hom(M(W),M(W))$ with $\Hom(A,A) \otimes \Hom(W,W)$ provides a map
\begin{equation}
                \phi\otimes \rho:V \otimes \Lg \to \Hom(M(W),M(W))
\end{equation}
which measures.  Moreover, $V\otimes \Lg$ is closed under Lie
bracket, as is $V_0 \otimes \Lg$.

        If $\kappa$ is the Killing form on $\Lg$ then $\tau \otimes \kappa$ is well defined on $V_0 \otimes \Lg$, and
$\tau \otimes \kappa([V_{\tau}\otimes \Lg, V_{0\tau}\otimes \Lg]) =
0$. Let $K(\Lg)$ be the kernel of $\tau \otimes \kappa$.  There is
then a short exact sequence of Lie algebras
\begin{equation}
        0 \to V_{0\tau} \otimes \Lg/K(\Lg) \to V_{\tau} \otimes \Lg/K(\Lg) \to V_{\tau}\otimes \Lg/V_{0\tau}\otimes \Lg \to 0
\end{equation}
The Lie algebra $V_{\tau} \otimes\Lg/K(\Lg)$ is the central
extension of the loop algebra
 $\Lg \otimes \Comp[x,x^{-1}] \approx V_{\tau} \otimes \Lg/V_{0\tau} \otimes\Lg$.
  It turns out that the cocycle $c$ of the central extension is
given by
\begin{equation}
        c(v \otimes \xi,w \otimes \psi) =  \tau\Omega (\mu v,\mu w)\kappa(\xi,\psi).
\end{equation}
 In the case of example i,  $V_{\tau} \otimes \Lg/V_{0\tau} \otimes \Lg = \Lg[x,x^{-1}]$,
  the loop algebra of $\Lg$, and $c$ is the expected central
extension,
  \begin{equation}
        c[x^m\xi,x^n\psi] = m \delta_{-m,n}\kappa(\xi,\psi).
  \end{equation}
 \end{enumerate}
 }
\end{Examples}
\section{Dual comodules, positive energy representations, and smooth
representations}

\subsection{Dual coalgebras and dual comodules}{\rm     If $A$ is an algebra, and
$M$ an $A$ module, the the constructions $P(A,\Comp), Q(M,\Comp)$,
have alternative descriptions which make hands on calculations
easy. }
  \begin{Proposition}\Hb {\rm
 \begin{enumerate}
  \renewcommand{\labelenumi}{(\roman{enumi})}
  \item $P(A,\Comp) =: A^{*} = \{\alpha: A\to \Comp: \Ker \alpha \geq I,
  I \mbox{ an ideal}, dim A/I< \infty\}$

   \item  $Q(M,\Comp) =: M^{*} = \{\mu: M \to \Comp: \Ker\mu \geq W, AW \leq W, dimM/W< \infty\}.$
 \end{enumerate}
  }\end{Proposition}
 \Proof.  i.Observe that since $A/I$ is finite dimensional
multiplication in $A/I$ gives the linear dual $(A/I)'$ the
structure of a coalgebra with the obvious measuring map into
$\Hom(A,\Comp)$. Then, since $(A/I)'$ maps to $P(A,\Comp)$ by the
universal property, $A^{*} = \lim(A/I)' \leq \Hom(A,\Comp)$ maps to
$P(A,\Comp)$.
        But now observe that the measuring map $\pi:P(A,\Comp) \to \Hom(A,\Comp)$ has
its image in $A^{*}$.  To see this consider $c$ in $P(A,\Comp)$.
Let $C$ be a finite dimensional subcoalgebra of $P(A,\Comp)$
containing $c$. Then the restriction of the measuring map $\pi:C\to
\Hom(A,\Comp)$ corresponds to an algebra homomorphism $\pi:A\to\Hom
(C,\Comp)$, this last being a finite dimensional algebra.  Let $J$
be the kernel of
 $\pi$. Since $\pi:A\to \Hom(C,\Comp)$ factors through $A/J$, $\pi(c):A \to \Comp$ must factor through $A/J$,
and $\pi(c)$ is in $A^{*}$ as required.

ii.  The argument is exactly parallel to that of (i).

 \begin{Remarks}\Hb{\rm
 \begin{enumerate}
  \renewcommand{\labelenumi}{(\roman{enumi})}
 \item Evidently $M^{*}$ becomes a module for the opposite algebra $A^{{\rm op}}$
under the action
\begin{equation}
                am = \sum_{(m)}m_{(1)}(a)m_{(0)}.
\end{equation}
One can ask what representations arise as dual comodules.  It is
evident that such a representation $V$ must have the property that
every element of $V$ lives in some finite dimensional submodule of
$V$. Representations which have this property will be called
locally finite.

 \item  More generally, given modules $M$, $N$ over $A$, $B$ respectively,
$Q(M,N)$ can be considered an $A^{{\rm op}}$ module.
\end{enumerate}

        The ingredients for the applications of interest are an algebra $A$ and a
representation of $A$ on a vector space $V$, and a distinguished
subalgebra $B$. Considered as an $A$ module, $V^{*} = ({}^AV)^{*}$
is not very interesting, and may in fact be zero.  However,
considered as a $B$ module, $({}^BV)^{*}$ is not only a $B$ module,
but an $A$ module. The property of the subalgebra $B$ which gives
$({}^BV)^{*}$ the structure of an $A$ module is as follows. }
\end{Remarks}
\begin{Definition} {\rm
 Say $B \leq A$ is {\em quasi-normal} if and only if for every $a$
in $A$ there exists $a_1,..., a_n$ such that
\begin{equation}
                BaB = \sum_1^{l} Ba_i  =  \sum_{1}^{l}a_iB.
\end{equation}
 }\end{Definition}
 \begin{Lemma}\label{QNAlem}  {\rm Suppose $B$ is quasi normal in $A$, and let $s:A \to A$
be an anti-automorphism.  Then if $M$ is a representation of $A$,
$Q({}^BM, \Comp)$ is an $A$ module.

 \Proof.    Notice that the action of $A^{{\rm op}}$ on $\Hom(M,\Comp)$ given
by
\begin{equation}
                a\mu(m) = \mu(am)
\end{equation}
coincides with the action of $B^{{\rm op}}$ on $({}^BM)^{*}$
whenever $\mu$ is in $({}^BM)^{*}$ and $a$ is in $B$.  Prefacing
this action with the antiautomorphism $s$
\begin{equation}
                a\Circ \mu(m) = \mu(s(a)m)
\end{equation}
defines an action of $A$ on \Hom(M,\Comp).  The claim is that
$({}^BM)^{*}$ is fixed by this action.

        Let $\alpha$ be in $({}^BM)^{*}$ and let $a$ be in $A$.
        By  4.2 $\alpha: M \to \Comp$ vanishes on
$N$, a $B$ submodule with $M/N$ finite dimensional.  The problem is
to show that there is a $B$ submodule $N_a$ such that $M/N_a$ is
finite dimensional and $(a\Circ \alpha)(N_a) = 0$.

        Since $B$ is quasi normal write
\begin{equation}
                Bs(a)B = \sum_{1}^l a_iB  =  \sum_1^l Ba_i.
\end{equation}
Then define linear maps
\begin{equation}
        \alpha_{i}~: N \to M \to M/N,    \quad \alpha_{i}(n) = a_in +N
        \end{equation}
and let
\begin{equation}
                N_i = \Ker \alpha_{i}, \quad N_a  = \cap N_i.
\end{equation}
Now observe that $N_a$ is a $B$ submodule of $M$: consider $a_ibn$
for $b$ in $B$, $n$ in $N_a$. We  can write
\begin{equation}
 a_ib = \sum_1^l b_j a_j
\end{equation}
so that
\begin{equation}
  \alpha_i  (bn) = a_i bn + N = \sum_1^l b_j a_j n +N.
\end{equation}
Since $n$ is in $N_a$, $a_jn$ is in $N$ for all $j$, and hence so
is $b_ja_jn$. Thus $bn$ is in the kernel of $a_i~$ for all $i$.

Finally check that $N_a$ is contained in $ker a \Circ \alpha$. For
$n$ in $N_a$\ $a \Circ\alpha(n) = \alpha(s(a)n)$. But $s(a)$ is in
 $Bs(a)B$, so $s(a) = \sum_1^l b_ja_j$ and $s(a)n = \sum_1^l b_ja_j  n$.
  Since $n$ is in $\Ker a_i~$ for each $i$, $a_in$ is in $N$ for
each $i$, hence $s(a)n$ is in $N$ and $\alpha(s(a)n) = 0$ as
required.
  }\end{Lemma}
\subsection{Application to totally disconnected groups}

        Let $G$ be a totally disconnected group (see [3] for a survey of the
representation theory of these objects) with a given compact open
subgroup $K$, and let $M$ be a complex representation of $G$, hence
a representation of $\Comp G (= A)$ and $\Comp K (=B)$.

\begin{Lemma}\label{QNlem}  {\rm $\Comp K$ is quasi normal in $\Comp G$.

\Proof.  Let $g$ be in $G$.  The double coset $KgK$ is a finite union of
either right or left cosets of $K$ and the left coset
representatives $\{g_i\}$ may be chosen to be the same as the right
coset representatives.  Then
\begin{equation}
                        \Comp K g\Comp K = \sum_{1}^n\Comp Kg_j = \sum_1^n g_j\Comp K
\end{equation}
as required.
  }\end{Lemma}
\begin{Corollary}\label{QREPCG} {\rm $Q({}^{\Comp K} M,\Comp)$ is a representation of $\Comp G$ which is locally
finite as a representation of $\Comp K$.

\Proof.  All that is needed to meet the conditions of lemma \ref{QNAlem}
 is the choice of an appropriate antiautomorphism.  Clearly the
map $s(g) = g^{-1}$ is a suitable choice.

        The representation $Q({}^{\Comp K} M,\Comp)$ is almost, but not quite the smooth dual
of $M$.  The relationship can be described in coalgebraic terms.
  }\end{Corollary}
\begin{Definition} {\rm If $F$ is a subcoalgebra of $C$, and $D$ is a $C$
comodule, define the {\em restriction} of $D$ to $F$ to be
\begin{equation}
                        {}_F|D = \{d \in D: \Delta d \in F \otimes D\}.
\end{equation}
Thus $ {}_F|D$   is a $C$ sucomodule and an $F$ comodule.

        In particular, the coalgebra $P(\Comp K, \Comp) = (\Comp K)^*$ contains as an important
subcoalgebra the vector space with basis $K^{\wedge}$, the set of
group homomorphisms $\rho:K \to \Comp$.  The trivial homomorphism
 $\tau:K \to \Comp$ in particular is in $K^{\wedge}$.  Consider the subcomodule
\begin{equation}
                {}_{\Comp\tau}|({}^{\Comp K} M).
\end{equation}
  }\end{Definition}
 \begin{Proposition}\Hb {\rm
\begin{enumerate}
 \renewcommand{\labelenumi}{(\roman{enumi})}
 \item  If $K'$ is another compact open subgroup of $G$ then $({}^{\Comp K} M)^{*} =  ({}^{\Comp K'}M)^{*}$.

  \item If $K'\leq K$, and if $\tau, \tau'$ are the corresponding trivial homomorphisms,
then
\begin{equation}
                {}_{\Comp\tau}|({}^{\Comp K} M)^* \leq {}_{\Comp\tau'}|({}^{\Comp K'}M)^*.
\end{equation}

 \item The union $\cup  {}_{\Comp\tau}|({}^{\Comp K} M)^*$ over all compact open $K$ is the smooth
dual of $M$.
\end{enumerate}
\Proof.  The only statement which is not immediate is the first.
Suppose that $K'\leq K$.  The inclusion induces a map
 $(\Comp K)^{*} \to (\Comp K')^{*}$, and any $(\Comp K)^{*}$ comodule is automatically a
 $(\Comp K')^{*}$ comodule. Moreover, any $(\Comp K)^{*}$ comodule $D$ equipped
with a measuring map
 $\rho:D \to \Hom(M,C)$ is also a measuring comodule for $(\Comp K')^{*}$.  Thus
 $({}^{\Comp K} M)^{*} \to (\Comp K'M)^{*}$.

        Less obviously $(\Comp K'M)^{*} \to ({}^{\Comp K} M)^{*}$.  Let $\alpha:M \to C$ vanish on $N'$
which is a $\Comp K'$ submodule with $M/N'$ finite dimensional. The
aim is to show that there exists $N$, a $\Comp K$ submodule with
 $\alpha(N) = 0$ and $M/N$ finite dimensional.

        Write $K'KK' = \sqcup k_iK' = \sqcup K'k_i$.  Since $K$ and $K'$ are compact open
this is a finite union.  The argument now is the same as that which
established \ref{QNlem}.  Define maps
\begin{eqnarray}
                k_i: N' &\to& M \to M/N' \End
                k_in' &=& k_n'+N' \mbox{\ for\ } n' \in N'
\end{eqnarray}
and set
\begin{equation}
  N= \cap \Ker k_i.
\end{equation}
The arguments that (i) $N$ is a $CK$ module, (ii) $N$ is contained
in $\Ker \alpha$ and (iii) $M/N$ is finite dimensional follow the
pattern of \ref{QNAlem}
  }\end{Proposition}
\subsection{Application to loop algebras}
(See [5] for basic information on the subject.) Let $\Lg$ be a
finite-dimensional simple Lie algebra over $\Comp$ and let
$\Lg[x,x^{-1}]$ denote the loop algebra of $\Lg$ consisting of
Laurent polynomials in $x$ with coefficients in $\Lg$. A
representation $M$ of $\Lg$ is a projective representation with
cocycle $c$ if
\begin{equation}
              (\xi x^i)(\psi x^j) m = \left( (\psi x^j)(\xi x^i) \right)
              +   [\xi,\psi]x^{i+j}m + c(\xi x^i,\psi x^j)m
\end{equation}
for all $m$ in $M$.  The representation is said to be of level $k$
if it is projective with cocycle $c$ given by
\begin{equation}
                        c(\xi x^i,\psi x^j) = ik\kappa(\xi,\psi)\delta_{i,-j}
\end{equation}
where $\kappa(\, ,\, )$ is the Killing form on $\Lg$.

A projective representation of $\Lg[x,x^{-1}]$ corresponds to an
ordinary representation of the central extension $\Lg[x,x^{-1}]
\oplus \Comp c$ in the usual way.  Thus level $k$ representations are
representations in which $c$ acts as multiplication by $k$.  In
addition, there is an outer derivation $d$ of $\Lg[x,x^{-1}]$ given
by
\begin{equation}
                                d\xi x^i = i \xi x^i.
\end{equation}
 Form the Lie algebra $\Lg[x,x^{-1}] \oplus \Comp c \oplus \Comp d$,
 setting $[d,\xi x^i] = i \xi x^i$, $[d,c]=0$.  The algebras of
interest are universal enveloping algebras of this Lie algebra and
certain subalgebras. Write
\begin{eqnarray}
               U &=& U(\Lg[x,x^{-1}]\oplus \Comp c \oplus \Comp d),\End
               U_{\geq} &=& U(\Lg[x] \oplus \Comp c \oplus \Comp d) \End
                U_{\leq} &=& U(\Lg[x^{-1}] \oplus \Comp c \oplus \Comp d) \End
                U_> &=& U(\Lg[x]x) \End
                U_< &=& U(\Lg[x^{-1}]x^{-1})
\end{eqnarray}
The isomorphisms as vector spaces
\begin{equation}
        \Lg[x,x^{-1}]  =  \Lg[x^{-1}]x^{-1} \oplus \Lg[x]
        = \Lg[x^{-1}]x^{-1} \oplus \Lg \oplus \Lg[x]x
\end{equation}
induce isomorphisms of vector spaces
\begin{equation}
           U  = U_< \otimes U_{\geq}
           = U_< \otimes U(\Lg \oplus \Comp c \oplus \Comp d) \otimes U_>.
\end{equation}
The bracket with $d$ provides a $\Bbc{Z}$ grading (as vector
spaces) of all the universal enveloping algebras described here.
With respect to this grading
\begin{equation}
                U_< = \oplus_{n \leq 0} (U_<)_n
\end{equation}
where $(U_<)_n$ is the set of elements of degree $n$.  Each
$(U_<)_n$ is finite dimensional.  Hence the subspace
\begin{equation}
            (U_<)_{(n)} = \oplus _{0 \leq j \leq n} (U_<)_j
\end{equation}
is also finite dimensional and
\begin{equation}
                U_< =  \oplus_{n\leq 0} (U_<)_{(n)}.
\end{equation}
\begin{Lemma}\label{UQNlem} {\rm $U_{\geq}$ is quasi-normal in $U$.

\Proof.  This is essentially a consequence of the analogue of the Poincare-
Birkhoff-Witt theorem for universal enveloping algebras.  Observe that
\begin{eqnarray}
        (U_<)_{(n)}\otimes U_{\geq}   &=&   U_{\geq} \otimes(U_<)_{(n)} \End
        U_{\geq} \otimes(U_<)_{(n)}   &=&   (U_<)_{(n)}\otimes U_{\geq}.
\end{eqnarray}
The result follows since $a$ in $U$ is in some $(U_<)_{(n)}
\otimes U_{\geq}$. If $\{a_{i}\}$ is a basis for $(U_<)_{(n)}$  then
\begin{equation}
                        U_{\geq}aU_{\geq}  =\sum_{i}a_iU_{\geq}  = \sum_{i} U_{\geq}a_i
\end{equation}
as required.

        The antiautomorphism commonly used is that determined by the Lie
algebra antiautomorphism $s:\Lg[x,x^{-1}] \oplus \Comp c \oplus
\Comp d
\to
\Lg[x,x^{-1}]\oplus \Comp c \oplus \Comp d$,
\begin{equation}
                        s(\xi x^i) = -\xi x^{-i}, s(c) = c, s(d) = -d.
\end{equation}
  }\end{Lemma}
 \begin{Proposition}\Hb {\rm   If $M$ is a $U$ module then
 \begin{enumerate}
 \renewcommand{\labelenumi}{(\roman{enumi})}
        \item $({}^{U_{\geq}}M)^*$ is a level $k$ representation if $M$ is.
        \item $({}^{U_{\geq}}M)^*$ is locally finite as a $U_{\geq}$ module.
 \end{enumerate}
\Proof.  i. This is more or less a direct corollary of \ref{UQNlem}.  Calculate
 \begin{eqnarray}
       \lefteqn{ [(\xi x^i)_(\psi x^j)\alpha - (\psi x^j)_(\xi x^i)\alpha - [x,\psi]xi+j\alpha ](m)} \End
        &=&  \alpha[(s(\psi x^j)(s(\xi x^i) - s(\xi x^i)s(\psi x^j) - s([x,\psi]x^{i+j}))m] \End
        &=&  \alpha[((\psi x^{-j})(\xi x^{-i}) - (\xi x^{-i})(\psi x^{-j})  - ([\psi,x]x^{-i-j}))m] \End
        &=&  \alpha[c(\psi x^{-j},\xi x^{-i})m] \End
        &=&  c(s(\psi x^j),s(\xi x^i))\alpha(m) \End
        &=&  -c(s(\xi x^i),s(\psi x^j))\alpha(m) \End
        &=&  c(\xi x^i,\psi x^j)\alpha_(m)
\end{eqnarray}
since $c(\xi x^i,\psi x^j) = ik\delta_{i,{-j}}\kappa(x,\psi) =
{-j}k\delta_{{-j},i}\kappa(\psi,x)
= c(\psi x^{-j},\xi x^{{-i}})$.  This establishes i.

For ii observe that $s(U_{\geq}) = U_{\geq}$.  The result then
follows from 4.6.
  }\end{Proposition}
\begin{Definition} {\rm Say a representation $M$ of $U$ is positive
energy if $d$ acts diagonally with real eigenvalues and the
eigenvalues of d are bounded above.
 }\end{Definition}
        As with the category of smooth representations of totally disconnected
groups, so the category of positive energy representations of a
loop algebra admits the existence of a dual.  As the smooth dual of
a representation of a totally disconnected group can be identified
in terms of restricted comodules, so the dual positive energy
representation of a representation $M$ can be identified as an
appropriate restriction of $({}^{U_{\geq}}M)^{*}$.  It remains to
identify the appropriate subcoalgebra of $(U_{\geq})^{*}$.

        The universal enveloping algebra $U_{>}$ has an augmentation ideal
$U_{>}{}^{+} = \oplus_{n>1} (U_{>})_n$. This generates an ideal
$U_0$ of $U_{\geq}$
\begin{equation}
                U_0 = U_{\geq} (U_{>}{}^{+}).
\end{equation}
A short calculation shows that $U_0$ is in fact a two sided ideal.
Define
\begin{equation}
                P_0{}^N   = im (U_{\geq}/(U_0)^N)^{*} \to (U_{\geq})^{*}.
\end{equation}
Since $P_0{}^{N+j} \geq  P_0{}^N$, define
\begin{equation}
                P_0 = \cup P_0{}^N.
\end{equation}
 \vfil\eject
\begin{Proposition} \Hb {\rm
\begin{enumerate}
 \renewcommand{\labelenumi}{(\roman{enumi})}
  \item ${}_{P_0}|({}^{U_{\geq}}M)^*$ is a $U$ submodule of $({}^{U_{\geq}}M)^{*}$.

  \item If ${}_{P_0}|({}^{U_{\geq}}M)^{*}$ is generated by a finite set of eigenvectors for $d$ then it
  is positive energy.

\end{enumerate}
 }\end{Proposition}
\Proof i.  Check that for $z$ in $U$, $q$ in ${}_{P_0}|({}^{U_{\geq}}M)^*$, $z\Circ q$ is
in $P_0 | ({}^{U_{\geq}}M)^*$, or equivalently, for some $N$, any
$u$ in $s(U_0{}^N)$, $u\Circ z\Circ q = 0$. Using $P_0 = \cup
P_0{}^N$ it can be shown that
\begin{equation}
                {}_{P_0}|({}^{U_{\geq}}M)^* = {}_{\cup P_0^N}|({}^{U_{\geq}}M)^*.
\end{equation}
Suppose then that $q$ is in $P_0^{N'}|({}^{U_{\geq}}M)^*$ for some
$N'$, that is, $u\Circ q = 0$ for all $u$ in $s(U_0{}^{N'})$.
        If $i \geq 0$ and $z$ in $(U_{\geq})i$, then $u\Circ z\Circ q = 0$ since $z$ is in
        $ U_{\geq}$ and $s(U_0{}^{N'})$ is an
ideal of $U_{\geq}$.  If $i< 0$, observe that
\begin{equation}
                s(U_0{}^N)(U_{\geq})_i  \leq (U_{\geq})_{(i)}s(U_0{}^{N+i}).
\end{equation}
Thus for $z$ in $(U_{\geq})_i$ , $q$ in
 $P_){}^{N'}|({}^{U_{\geq}}M)^{*}$,  $u\Circ z \Circ q  = 0$ provided $N> N'{{-i}}$.

 ii.  Write $V = {}_{P_0}|({}^{U_{\geq}}M)^*$.   Assume that $\{q_i\}$ is a finite generating set
for $V$ of $d$ eigenvectors.  Since $V$ is locally finite as a
$U_{\geq}$ module,  $U_{\geq}\{q_i\}$ is a finite dimensional
$U_{\geq}$ module, call it $D$. In particular, the element $d$ acts
on $D$, and is diagonalizable on $D$ with finitely many
eigenvalues. But then since
\begin{equation}
                V = UD = U_{<}D
\end{equation}
$d$ acts diagonally on $V$ and the eigenvalues are bounded below.

\vskip 0.1cm
{\bf References}
\vskip 0.1cm
\begin{enumerate}
\item Batchelor, M. In search of the graded manifold of maps between
graded manifolds.  In {\em Complex Differential Geometry and
Supermanifolds in Strings and Fields.}  P.J.M. Bongaarts and R.
Martini, eds. LNP 311. Springer, 1988.

\item Batchelor, M. Measuring coalgebras quantum group-like objects,
and non-commutative geometry.  In {\em Differential Geometric
Methods in Theoretical Physics.} C. Bartocci, U. Bruzzo and R.
Cianci, eds. LNP 375. Springer. 1990.

\item Cartier, P. Representations of p-adic groups: a survey. In {
\em Automorphic forms and L-functions.  Proceedings of Symposia in Pure
Mathematics.}  Volume XXXIII, part 1. p 111-155. AMS. Providence,
1979.

\item Hyland, M. and Batchelor, M. Enrichment in coalgebras and
comodules: an approach to fibrations in the erniched context. {\em
In preparation.}

\item Kac, V. G.  {\em Infinite dimensional Lie algebras.}
Birkhauser, Boston, 1983.

\item Kelly, G. M.  {\em Basic Concepts of Enriched Category Theory.}
London Mathematical Lecture Note Series 64.  Cambridge University
Press, 1985.

\item Kobayashi, S. and Nomizu,K.  {\em Foundations of Differential
Geometry.} Wiley. 1996.

\item Spivak, M. {\em A Comprehensive Introduction to Differential
Geometry,} vol. 2. 2nd ed. Publish or Perish, 1979.

\item Sweedler, M. {\em Hopf Algebras.}  Benjamin. New York, 1969.
\end{enumerate}
\end{document}